\documentclass[a4paper]{amsart}
\usepackage{amssymb}
\usepackage[all]{xy}
\newtheorem{theorem}{Theorem}[section]

\newtheorem{lemma}[theorem]{Lemma} 

\newtheorem{corollary}[theorem]{Corollary}
\theoremstyle{definition}

\newtheorem{definition}{Definition}[section]
\theoremstyle{remark}

\newtheorem{remark}[theorem]{Remark}

\newtheorem{example}{Example}
\newlength\xhomlength
\DeclareMathOperator\Forb{{\rm Forb}_h}
\newcommand{\field}[1]{\mathbb{#1}}
\newcommand{\bbbn}{\field{N}}

\newcommand{\xhom}[1]{
  \settowidth{\xhomlength}{$#1$}
\ \xy\xymatrix@+=\xhomlength{*!{}\ar@{->}[r]^{#1}&*!{}}\endxy\ }
\newcommand{\nxhom}{\
  \xy\xymatrix{*!{}\ar@{->}[r]\ar@{{}|}[r(0.45)]&*!{}}\endxy\ }
\newcommand{\vxhom}[1]{
  \settowidth{\xhomlength}{$#1$}
  \ \xy\xymatrix@+=\xhomlength{*!{}\ar@{.>}[r]^{#1}&*!{}}\endxy\ }

\DeclareMathOperator\depth{{\rm td}}
\DeclareMathOperator\clos{{\rm clos}}
\DeclareMathOperator\height{{\rm height}}
\DeclareMathOperator\tw{{\rm tw}}

\DeclareMathOperator\mad{{\rm mad}}

\newcommand{\rdens}[1]{\nabla_{#1}}
\newcommand{\dens}{\nabla}
\newcommand{\card}[1]{\lvert{#1}\rvert}

\newcommand{\rdpath}[1]{\xy\xymatrix{*!{}\ar@{~2>}[r]^{#1}&*!{}}\endxy}
\newcommand{\ldpath}[1]{\xy\xymatrix{*!{}&*!{}\ar@{~2>}[l]_{#1}}\endxy}

\begin{document}
\title[Bounded Expansion III: Homomorphism Dualities]{Grad and Classes 
with Bounded Expansion III.\\
Restricted Graph Homomorphism Dualities} \author{Jaroslav Ne\v
  set\v ril}
\address{Department of Applied Mathematics\\
  and\\
  Institute of Theoretical Computer Science (ITI)\\
  Charles University\\
  Malostransk\' e n\' am.25, 11800 Praha 1\\
  Czech Republic} \email{nesetril@kam.ms.mff.cuni.cz} \author{Patrice
  Ossona de Mendez}
\address{Centre d'Analyse et de Math\'ematiques Sociales\\
  CNRS, UMR 8557\\
  54 Bd Raspail, 75006 Paris\\
  France} \email{pom@ehess.fr}
\begin{abstract}
We study restricted homomorphism dualities in the context of classes with bounded expansion.
This presents a generalization of restricted dualities obtained earlier for bounded degree graphs and also for proper minor closed classes.
This is related to distance coloring of graphs and to the ``approximative 
version'' of Hadwiger conjecture.
\end{abstract}
\maketitle
\section{Introduction}

We motivate this paper by the following two examples.

$\bf{Example 1.}$
Celebrated Gr\" otzsch's theorem (see e.g. \cite{D}) says that every 
planar graph is $3$-colourable. In the language of homomorphisms this 
says that for every triangle free planar graph $G$ there is a homomorphism 
of $G$ into $K_3$. Here
a {\em homorphism} from a graph $G$ to a graph $H$ is a mapping
$f:V(G)\rightarrow V(H)$ which preserves adjacency: $ \{f(x),f(y)\}\in E(H)$ 
whenever $\{x,y\}\in
E(G)$.  $G\xhom{f}H$ or $f: G\xhom{}H$ denotes that $f$ is a homomorphism from $G$ to $H$. The existence of a homomorphism from $G$
to $H$ is noted $G\xhom{} H$, while the non-existence of such a
homomorphism is noted  $G\nxhom H$. It is also clear that the relation 
$G \leq H$ defined as $G\xhom{} H$ is a quasiorder on the class of all finite graphs. This quasiorder
becames partial order if we restrict it to the class of all minimal retracts (i.e. {\em cores}). This partial order is called {\em homomorphism order}.
See \cite{HN} for a recent introduction to graphs and homomorphisms.

Using the partial order terminology the Gr\" otzsch's theorem 
says that $K_3$ is an upper bound (in the homomorphism order) for 
the class  $\mathcal P_3$ of all planar triangle free  graphs. As obviously 
$K_3 \not\in \mathcal P_3$ a natural question (first formulated in \cite{N1})
suggests: Is there yet a smaller bound? The answer, which may be viewed as 
a strengthening of Gr\" otzsch's theorem, is positive. Thus there exists 
a triangle free $3$-colorable  graph $H$ such that $G \xhom{} H$ for every graph $G\in \mathcal P_3$. This has been proved in 
\cite{Taxi_tdepth,fold} in a stronger version for minor closed classes. 
The case of planar graphs and triangle is interesting in its own 
and it has been related to the Seymour conjecture and Guenin's theorem 
\cite{Gue}, see \cite{Nas} and seems to  found a proper setting in the context of 
$TT$-continuous mappings, see \cite{NSam}.
Restricted duality results have been generalized since to proper 
minor closed  classes of graphs and to other forbidded subgraphs. In fact to 
any finite set of connected graphs, see \cite{Taxi_tdepth}.
This then implies that Gr\" otzsch's theorem can be strengthened by a
sequence of ever stronger bounds and that the supremum of the class of
all triangle free planar graphs does not exist, see \cite{bound}. 
  
$\bf{Example 2.}$ 
Let us consider all sub-cubic graphs (i.e. graph with maximum degree $\leq 3$).
By Brooks theorem (see e.g. \cite{D}) all these graphs are $3$-colorable 
with the single connected exception $K_4$. What about the class of all 
sub-cubic {\em triangle free} graphs? Does there exists a triangle free 
$3$-colorable bound? The positive answer to this question is given in \cite{DMN} and \cite{HH}.
In fact for every finite set ${\mathcal F} = \{F_1, F_2, \ldots, F_t\}$ of 
connected graphs
there exists a graph $H$ with the following properties:

- $H$ is $3$-chromatic;

- $G \xhom{} H$ for every subcubic graph $G \in \Forb(\mathcal F)$.

(Here $\Forb(\mathcal F)$ is the class of all  graphs $G$ which satisfy
$F_i \nxhom{} G$ for every $i = 1,2,\ldots, t$.)
In this case we briefly say that the class of all sub-cubic graphs has all
{\em restricted dualities}.  (We shall motivate this terminology below.)

It is interesting to note that while sub-cubic graphs have restricted 
dualities 
(and, more generally, this also holds for the classes of bounded degree graphs)
for the classes of degenerated graphs a similar statement is not true
(in fact, with a few trivial exceptions, it is never true), 
see \cite{N1,bound}.

Where lies the boundary for validity of restricted dualities?
This is the central question of this paper. We give a very general 
sufficient condition for a class to have restricted dualities. 
But first we introduce another source for restricted dualities. Chronologically this is also the original cotext.
 
The following is a  partial order formulation of an  
important homomorphism (or coloring) problem:

\begin{definition}
A pair $F,D$ of graphs is called {\em dual pair} if for every graph $G$ holds:
\begin{equation}
\label{eq1}
F \nxhom{} G  \iff  G \xhom{} D.
\end{equation}

\end{definition}
We also say that $F$ and $D$ form a duality, $D$ is called {\em dual of} $F$.
Dual pairs of graphs and even  of relational structures were characterized
in \cite{NT}, the notion itself goes back to \cite{NP}.
Equivalently, one can describe a dual pair $F,D$ by saying that for the class
$\Forb(F)$ the graph $D$ is the maximum graph (in the homomorphism order).
 
It appears (and this is the main result of \cite{NT}) that
(up to the homomorphism equivalence) 
all the dualities are of the form $(T, D_T)$  where 
$T$ is a finite (relational) tree. Every  dual $D_T$ is  uniquelly determined by the  tree $T$ 
(but its structure is by far more difficult to describe,
see e.g. \cite{NT1,NS}).
These results imply in most cases infinitely many examples. But  a 
much richer spectrum (and in fact a surprising richness of results) is 
obtained by restricting the validity of  \eqref{eq1}
to a particular class of graphs $\mathcal K$:  

\begin{definition}
A class $\mathcal K$ admits a {\em restricted duality} if, for any
finite set of connected  graphs $\mathcal F = \{F_1,F_2,\ldots,F_t\}$, 
there exists a finite graph $D_{\mathcal F}^{\mathcal K}$ such
that $F_i\nxhom D_F^{\mathcal K}$ for $i=1,\ldots,t$ and such that
for all $G\in\mathcal K$ holds:
$$(F_i\nxhom G), i = 1,2,\ldots,t, \iff (G\xhom{} D_F^{\mathcal K})$$.
\end{definition}

It is easy to see that using the  homomorphism order we can reformulate this definition as follows: A class $\mathcal K$ has restricted duality if for any
finite set of connected  graphs $\mathcal F = \{F_1,F_2,\ldots,F_t\}$
the class $\Forb(\mathcal F) \cap \mathcal K$ 
is bounded in the class $\Forb(\mathcal F)$.
 
In our companion papers \cite{POMNI,POMNII} we defined the notion of grad 
and bounded expansion class. For the benefit of the reader we recall these
definitions in Section 2. The following is then the main result of this 
paper:

\begin{theorem}
\label{main}
Any class of graphs with bounded expansion has all restricted dualities.
\end{theorem}

As both proper minor closed classes and bounded degree graphs form 
classes of bounded expansion this result generalizes both Examples 1. and 2.
In fact the seeming incomparability of bounded degree 
graphs and minor closed classes led us to the definition of bounded 
expansion classes.

This paper is organized as follows. In Section 2 we recall basic definitions 
and results of \cite{POMNI} which will be needed.
In Section 3 we reformulate the restricted dualities in terms of local 
homomorphism properties and introduce the basic construction. 
In Section 4 we prove Theorem \ref{main} and in 
Section 5 we list several corollaries. Among them is a surprising result 
that exact odd powers of graphs in any given bounded expansion class 
have bounded chromatic number.

\section{Bounded expansion classes.}
In \cite{Taxi_tdepth}, we introduced
the {\em tree-depth} $\depth(G)$ of a graph $G$ as follows:

A {\em rooted forest} is a disjoint union of rooted trees.
The {\em height} of a vertex $x$ in a rooted
forest $F$ is the number of vertices of a
path from the root (of the tree to which $x$ belongs to) to $x$ and is noted $\height(x,F)$.
The {\em height} of $F$ is the maximum height of the vertices of $F$.
Let $x,y$ be vertices of $F$. The vertex $x$ is an {\em ancestor} of $y$ in $F$ if $x$ belongs to
the path linking $y$ and the root of the tree of $F$ to which $y$ belongs to.
The {\em closure} $\clos(F)$ of a rooted forest $F$ is the graph with
vertex set $V(F)$ and edge set $\{\{x,y\}: x\text{ is an ancestor of
  }y\text{ in }F, x \neq y\}$. A rooted forest $F$ defines a partial order on its set of vertices:
$x\leq_F y$ if $x$ is an ancestor of $y$ in $F$.
The comparability graph of this partial order is obviously $\clos(F)$.
The {\em tree-depth} $\depth(G)$ of a graph $G$ is the minimum height
of a rooted forest $F$ such that $G\subseteq\clos(F)$. As a
consequence, we have an algorithmic definition of the tree depth :

\begin{lemma}[\cite{Taxi_tdepth}]
\label{lem:recur}
Let $G$ be a graph and let $G_1,\dots,G_p$ be its connected components. Then:
\begin{equation*}
\depth(G)=\begin{cases}
1,&\text{if }|V(G)|=1;\\
1+\min_{v\in V(G)}\depth(G-v),&\text{if }p=1\text{ and }|V(G)|>1;\\
\max_{i=1}^p \depth(G_i),&\text{otherwise.}
\end{cases}
\end{equation*}
\end{lemma}

We say that
a class $\mathcal C$ has a {\em low tree-depth coloring} if, for any
integer $p\geq 1$, there exists an interger $N(p)$ such that any graph
$G\in\mathcal C$ may be vertex-colored using $N(p)$ colors so that
each of the connected components of the subgraph induced
by any $i\leq p$ parts has tree-depth at most $i$.
As it obviously holds $\depth(G)\geq\tw(G)-1$ any class having a 
low-tree depth coloring
has also low tree-width coloring (in the sense of \cite{2tw}.

The existence of low-tree depth colorings is related to the ntion of $p$-centered coloring, 
which have also been
introduced in \cite{Taxi_tdepth}:
A {\em $p$-centered coloring} of a graph $G$ is a vertex coloring
such that, for any  connected subgraph $H$, either some color $c(H)$
appears exactly once in $H$, or $H$ gets at least $p$ colors.
For the sake of completeness we recall some results of \cite{Taxi_tdepth}. 
These statements  establish the relationship of centered colorings 
and low tree/depth colorings. They are easy to prove (with the exception of 
Theorem \ref{mainold1} which is the central result of \cite{Taxi_tdepth}):

\begin{lemma} [\cite{Taxi_tdepth}]
\label{lem:centerlocal}
Let $G, H$ be graphs, let $p=\depth(H)$, let $c$ be a
$q$-centered coloring of $G$ where $q\geq p$.
Then any subgraph $H'$ of $G$ isomorphic to $H$ gets at least $p$
colors in the coloring $c$ of $G$.
\qed
\end{lemma}

From this lemma follows that $p$-centered colorings induce low
tree-depth colorings:

\begin{corollary}
\label{cor:c2t}
Let $p$ be an integer, let $G$ be a graph and let $c$ be a
$p$-centered coloring of $G$.

Then $i<p$ parts induce a subgraph of tree-depth at most $i$
\end{corollary}
\begin{proof}
Let $G'$ be any subgraph of $G$ induced by $i<p$ parts. 
Assume $\depth(G')>i$. According to Lemma \ref{lem:recur}, the
deletion of one vertex decreases the tree-depth by at most one. Hence
there exists an induced subgraph $H$ of $G'$ such that
$\depth(H)=i+1\leq p$.
According to lemma \ref{lem:centerlocal}, $H$ gets
at least $p$ colors, a contradiction.
\end{proof}

\begin{theorem}
\label{mainold1}
Any graph $G$ has $p$-centered coloring for any $p \leq td(G)$.
\end{theorem}

The following was established  in \cite{Taxi_tdepth} for the case of proper
minor closed classes of graphs. We prove it here in full generality.

\begin{theorem}
\label{th:col}
Let $\mathcal C$ be a class of graphs having low tree-depth colorings
and let $p$ be an integer. Then there exists integer $X(p)$,
such that every graph in $\mathcal C$ has a $p$-centered
coloring using $X(p)$ colors.
\end{theorem}

\begin{proof}
Let $G\in\mathcal C$. According to the assumption, there exists a
vertex partition into 
$C(p)$ parts, such that any $p$ parts form a graph of
tree-depth at most $p$. This partition will be defined as a
coloring $\bar{c}: V(G) \longrightarrow \{1,2,\ldots,C(p)\}$. For
any set $P$ of $p$ parts let $G_P$ be the graph induced by all the
parts in $P$. It is easy to see that any graph $G$ with $\depth(G) \leq p$
has a $p$-centered coloring $c$ by $\depth(G($ colors: 
we simply assign to any connected subgraph $H$ of $G$
the minimal level of a vertex of $H$ in the tree $F$ satisfying $H \subset \clos(F)$
(see the definition of tree depth at the beginning of this section.  
Consider the following (``product'') coloring $c$ defined as
\[
c(v) = (\bar{c}(v), (c_P(v); |P| = p, P \subset \{1,2,\ldots,C(p)\})).
\]

Take the product of the coloring of $G$ by $C(p)$ colors and of the
colorings of the $G_P$ as a new coloring of $G$ (with $X(p)=C(p)
N(p,p-1)^{\binom{C(p)}{p}}$ colors). Let $H$ be a connected
subgraph of $G$. Then, either $H$ gets at least $p+1$ colors, or
$V(H)$ is included in some subgraph $G_P$ of $G$ induced by $p$
parts. In the later case, some color appears exactly once in $H$.
\end{proof}

Recall that the {\em maximum average degree} $\mad(G)$ of a graph $G$ is
the maximum over all subgraphs $H$ of $G$ of the average degree of $H$, that is
$\mad(G)=\max_{H\subseteq G}\frac{2|E(H)|}{|V(H)|}$. The {\em
  distance} $d(x,y)$ between two vertices $x$ and $y$ of a graph is the minimum
length of a pth linking $x$ and $y$, or $\infty$ if $x$ and $y$ do not 
belong to  same connected component.

We introduce several notations:

\begin{itemize}
\item 
The {\em radius} $\rho(G)$ of a connected graph $G$ is:
\begin{equation*}
\rho(G)=\min_{r\in V(G)}\max_{x\in V(G)} {\rm d}(r,x)
\end{equation*}
\item
 A {\em center} of $G$ is a vertex $r$ such that $\max_{x\in V(G)}{\rm
  d}(r,x)=\rho(G)$.
\end{itemize}

\begin{definition}
Let $G$ be a graph. A {\em ball} of $G$ is a subset of vertices
inducing a connected subgraph. 
The set of all the families of pairwise disjoint balls of $G$ is noted $\mathfrak{B}(G)$.

Let $\mathcal P=\{V_1,\dotsc,V_p\}$ be a family of pairwise disjoint 
balls of $G$.
\begin{itemize}
\item
 The {\em radius} $\rho(\mathcal P)$ of $\mathcal P$ is
 $\rho(\mathcal P)=\max_{X\in \mathcal P}\rho(G[X])$ 
\item
The {\em quotient} $G/\mathcal P$ of $G$ by 
  $\mathcal P$ is a graph with vertex set $\{1,\dotsc,p\}$ and edge
  set $E(G/\mathcal P)=\{\{i,j\}: (V_i\times V_j)\cap
  E(G)\neq\emptyset\text{ or }V_i\cap V_j\neq\emptyset\}$.
\end{itemize}
\end{definition}

We introduce several invariants that generalize the one of maximum
average degree:
\begin{definition} The {\em greatest reduced average density (grad)}of $G$
      with rank $r$ is  $$\rdens{r}(G)=
\max_{\substack{\mathcal P\in\mathfrak{B}(G)\\
\rho(\mathcal P)\leq r}}\frac{|E(G/\mathcal P)|}{|\mathcal P|}$$ 

For the sake of simplicity, we also define:

The {\em grad of $G$}: $$\dens(G)=\max_{r}\rdens{r}(G)=\max_{H\preceq
    G}\frac{|E(H)|}{|V(H)|}$$

\end{definition}

Notice the two following well known facts (usually expressed by mean
of the maximum average degree):
\begin{lemma}
\label{fact:md}
Let $G$ be a graph. Then $G$ has an orientation such that the maximum
indegree of $G$ is at most $k$ if and only if $k\geq \rdens{0}(G)$.
\end{lemma}
\begin{lemma}
Let $G$ be a graph. Then $G$ is $\lfloor
2\rdens{0}(G)\rfloor$-degenerated, hence 
$\lfloor 2\rdens{0}(G)+1\rfloor$-colorable.
\end{lemma}

The following is our key definition:
\begin{definition}
\label{maindef}
A class of graphs $\mathcal C$ has {\em bounded expansion} if there
exists a function $f:\bbbn\rightarrow\bbbn$ such that for every graph 
$G \in\mathcal C$ and every $r$ holds
\begin{equation}
 \rdens{r}(G)\leq f(r).
\end{equation}

$f$ is called the expansion function.
\end{definition}

The following is a special case of the main
result  of        \cite{POMNI}
\begin{theorem}
\label{mainold}
For a class $\mathcal C$ of graphs are the following statements equivalent

\begin{itemize}
\item $\mathcal C$ has bounded expansion,
\item $\mathcal C$ has low tree-depth colorings,
\item $\mathcal C$ has $p$-centered coloring for every $p \geq 1$.       
\end{itemize}

\end{theorem}

\section{A construction}

\begin{definition}
Let $G,H$ be graphs and let $\mathcal P$ be a system 
of subsets of $V(G)$. We say that
$G$ is {\em $\mathcal P$-locally homomorphic to $H$} and 
we denote $G\vxhom{\mathcal P}H$ if for every 
subset $A \in \mathcal P$:
$$G[A]\xhom{} H$$.

We shall deal mostly with the following systems:
If $\phi: V(G) \longrightarrow X$ is a function and $p$ a positive integer
then we can consider the system $\mathcal P = \{A \subset V(G); |\phi(A)| \leq p\}$. 
This system will be denoted by $\mathcal P_{(\phi,p)}$. In this case 
we also say that $G$ is $(\phi,p)$-locally homomorphic to $H$ 
(instead of  $\mathcal P_{(\phi,p)}$-locally homomorphic). This is also denoted by 
$G \vxhom{(\alpha, p)} H$.
\end{definition}

\begin{example}
For  $H = K_2$ and $\phi$ an identical map $V(G) \longrightarrow V(G)$ 
a graph $G$ is $(\phi,p)$-locally homomorphic to $H$ iff the odd-girth of $G$
is $> p$.
\end{example}

The following is a modification of a construction introduced in \cite{fold}:  

\begin{definition}
Let $G,H$ be finite graphs and let $1\leq p<\card{V(H)}$ be an integer.
For $v\in V(H)$, define $\mathcal A_v=\{(I, v); I \in\binom{V(H)}{p}: v \in
I\}$, where $\binom{V(H)}{p}$ stands for the subsets of
$V(H)$ with cardinality $p$.
Define the sets $V_v=V(G)^{\mathcal A_v},
W=\bigcup_{v\in V(H)}V_v$ and the function $\alpha:W\rightarrow
V(H)$ by $\alpha(z)=v$ if $z\in V_v$.

The {\em $p$-truncated $H$-power} $G^{\Uparrow_p^H}$ of $G$ is the
graph with vertex set $W$ and with the edge set $F$ defined as follows:
$\{z,z'\} \in F$ iff $z \in V_v, z' \in V_{v'}$ and for every 
$I \in\mathcal A_v\cap\mathcal A_{v'}$ holds $\{z_{(I,v)},z_{(I,v')}'\}\in
E(G)$.
\end{definition}

\begin{remark}
$K_1^{\Uparrow_1^H}$ is isomorphic to $H$.
\end{remark}

\begin{remark}
The order of $G^{\Uparrow_p^H}$ is
$\card{V(H)}.\card{V(G)}^{\binom{\card{V(H)}-1}{p-1}}$.
\end{remark}

The function $\alpha:V(G^{\Uparrow_p^H})\rightarrow V(H)$ is
called the {\em color projection} of $G^{\Uparrow_p^H}$. This is justified by
\begin{lemma}
The color projection $\alpha$ of $G^{\Uparrow_p^H}$ is a homomorphism
from $G^{\Uparrow_p^H}$ to $H$:
$$\alpha: G^{\Uparrow_p^H}\xhom{}H$$
\end{lemma}
\begin{proof}
By definition, for any edge $\{z,z'\}$ of $G^{\Uparrow_p^H}$, there
exists an edge $\{v,v'\}$ of $H$ such that $z\in\mathcal A_v$ and
$z'\in\mathcal A_{v'}$, i.e. $v=\alpha(z)$ and $v'=\alpha(z')$.
\end{proof}

\begin{lemma}
Let $\alpha$ be the color projection of $G^{\Uparrow_p^H}$. Then
the graph $G^{\Uparrow_p^H}$ is $(\alpha,p)$-locally homomorphic to $G$:
$$G^{\Uparrow_p^H}\vxhom{(\alpha,p)}G$$. 
\end{lemma}
\begin{proof}
Let $A$ be a subset of $V(G^{\Uparrow_p^H})$ such that
$\card{\alpha(A)}\leq p$. Let $I$ 
be any subset of $V(H)$ of
cardinality $p$ such that  $\alpha(A) \subset I$. According to the definition
of $G^{\Uparrow_p^H}$, $\{z_{(I,\alpha(z))},z'_{(I,\alpha(z'))}\}\in E(G)$ for any $\{z,z'\}\in
E(G^{\Uparrow_p^H}[A])$. It follows that the mapping $z\mapsto z_{(I,\alpha(z))}$ is
a homomorphism from $G^{\Uparrow_p^H}[A]$ to $G$.
\end{proof}
 
\begin{lemma}
\label{lem:power}
Let $G,H,U$ be finite graphs, let $p$ be an integer.
Assume that $\gamma: G \xhom{} H$ is a homomorphims and that $G$ is $(\gamma,p)$-locally homomorphic to $U$. Schematically:

%neeed change!!!

Assume $\xy\xymatrix{G\ar[r]^\gamma\ar@{.>}[d]_{(\gamma,p)}&H\\
U}\endxy$

Then there exists a homomorphism $f: G \xhom{}U^{\Uparrow_p^H}$ 
such that $\gamma=f\circ\alpha$. Moreover, $U^{\Uparrow_p^H}$ is 
$(\alpha,p)$-locally homomorphic to $U$.

Schematically, this can be expressed by the following scheme:
%need change!!!!

$$\xy\xymatrix@C+0.5cm{
G\ar[r]^\gamma\ar@{.>}[d]_{(\gamma,p)}\ar[dr]^f
&H\\
U&U^{\Uparrow_p^H}\ar[u]_{\alpha}
\ar@{.>}[l]^{(\alpha,p)}}\endxy$$
\end{lemma}
\begin{proof}
For $I\in\binom{V(H)}{p}$ put  $G_I=G[\gamma^{-1}(I)]$ and let  $g_I$ be a
homomorphism from $G_I$ to $U$. Define $f$ as follows:
Given $x\in V(G)$ we define $f(x) \in gamma(x) $  by the following formula
$f(x)(I,\gamma(x)) = g_I(x)$ (see the above definition of $\Uparrow_p^H$).
Obviously  $f\circ\alpha = \gamma$. We prove that $f$ is a homomorphism.

Let $\{x,y\}$ be any edge. It is 
$\{\gamma(x),\gamma(y)\}\in E(H)$ as
$\gamma:
G\xhom{}H$. For any $I$ which contains both $\gamma(x)$ and $\gamma(y) $ holds
$\{f(x)_{(I,\gamma(x))},f(y)_{(I,\gamma(y)}=\{g_I(x),g_I(y)\}\in E(U)$. 
It follows that
$\{f(x),f(y)\}$ is an edge of $U^{\Uparrow_p^H}$ and thus
 $f$ is a
homomorphism. 
\end{proof}
This lemma  highlights a fundamental property of
$G^{\Uparrow_p^H}$ which we will state as follows:
\begin{lemma}
\label{lem:locbound}
Let $G,H,U$ be finite graphs and let $p$ be an integer. Then:
there is a homomorphism $G \xhom{} U^{\Uparrow_p^H}$ iff
there exists a homomorphism 
$\gamma: G \xhom{} H$ and $G$ is $(\gamma,p)$-locally homomorphic to $U$. 
Schematically, 
this may be depicted as follows:
$$G\xhom{}U^{\Uparrow_p^H} \quad\iff\quad
\xy\xymatrix{
G\ar[r]^\gamma\ar@{.>}[d]_{(\gamma,p)}&H\\U}\endxy$$
\end{lemma}
\begin{proof}
First,
assume $f:
G\xhom{ }U^{\Uparrow_p^H}$. Let $\alpha$ be the color
projection of $U^{\Uparrow_p^H}$ to $H$. Put $\gamma=\alpha \circ f$.
We have $\gamma:
G\xhom{}H$. Let $A\subseteq V(G)$. The condition
$\card{\gamma(A)}\leq p$ is equivalent to the condition
$\card{\alpha(f(A))}\leq p$. Hence the homomorphism
$f:
G\xhom{}U^{\Uparrow_p^H}$ together with the $(\alpha,p)$-local homomorphism 
of $\Uparrow_p^H$ to $U$ implies $(\gamma,p)$-local homomorphism 
of $\Uparrow_p^H$ to $U$. The reverse implications follows from the previous lemma.
\end{proof}

It is interesting to note that if we consider $G=U$ we get:
\begin{corollary}
$ G\xhom{}G^{\Uparrow_p^H}\iff G\xhom{}H$. In particular, $G$ is
homomorphism-equivalent to $G^{\Uparrow_p^G}$.
\end{corollary}

\begin{theorem}
\label{locdual}
A class of graphs $\mathcal C$ has restricted dualities 
iff
for any finite set $\mathcal F$ of graphs
there exist a graph $H$ and a graph $U\in\Forb(\mathcal F)$ such that
for every  $G\in\mathcal C\cap\Forb(\mathcal F)$ there exists a homomorphism 
$\gamma: G \xhom{} H$ for which $G$ is $(\gamma,p)$-locally 
homomorphic to $U$, where $p = \max\{\card{V(F)}; F \in {\mathcal F}\}$.

\end{theorem}

\begin{proof}
If $\mathcal C$ has restricted dualities and a set  $\mathcal F$ of graphs 
has a
dual $D_{\mathcal F}^{\mathcal C}$ then we may put  $U=H=D_{\mathcal F}^{\mathcal C}$.

Now assume that graphs  $U, H$ and a homomorphism $\gamma$ exist. Put  
$p=\max \{ \card{V(F)}; F \in \mathcal F\}$. 
In this situation we  prove that 
$U^{\Uparrow_p^{H}}$ is a dual of $\mathcal F$.
Then  for every 
$F \in \mathcal F$ holds $F\nxhom U^{\Uparrow_p^{H}}$. (Suppose contrary, let 
$F\xhom{g} U^{\Uparrow_p^{H}}$. By Lemma \ref{lem:power} $U^{\Uparrow_p^H}$ is $(p,p)$-locally homomorphic to $U$ and this together with 
$\card{g(V(F))}\leq\card{V(F)} \leq p$ would imply $F\xhom{} U$.
If $F\nxhom G$ for every $F \in \mathcal F$ then
$G\xhom{} U^{\Uparrow_p^{H}}$ according to Lemma~\ref{lem:locbound}. 
If $G\xhom{} U^{\Uparrow_p^{H}}$ then $F\nxhom G$ as
$F\xhom{} G$ would imply $F\xhom{} U^{\Uparrow_p^{H}}$.
\end{proof}

\section{Restricted dualities}

We shall need one more (``finitness'') result proved in \cite{Taxi_tdepth}, Corollary 3.3:
\begin{lemma}
\label{finitness}
For any positive integer $p$ there exists a number $\digamma(G))$  any graph $G$ with 
$\depth(G) \leq p$ is hom-equivalent to one of 
its induced
subgraph of order at most $\digamma(\depth(G))$.
\end{lemma}

\begin{theorem}
\label{main1}
Let $\mathcal F$ be a finite set of finite connected  graphs.
Then, for any  class  of graph $\mathcal K$ with bounded expansion
there exists a finite graph $U(\mathcal K,\mathcal F)\in\Forb(\mathcal
F)$ such that any graph of $\mathcal K\cap\Forb(\mathcal F)$ has
a homomorphism to $U(\mathcal K,\mathcal F)$.
\end{theorem}
\begin{proof}
Let $p=\max_{F\in\mathcal F}|V(F)|+1$. There exists an integer $N$, such that any graph
$G\in\mathcal K$ has a proper $N$-coloring in which any $p$ colors
induce a graph of tree depth at most $p$.
According Lemma~\ref{finitness}, there exists a finite set
$\hat{\mathcal D}_k$ of graphs with tree depth at most $k$, so that
any graph with tree-depth at most $k$ is hom-equivalent to one graph
in the set. Let $U(\mathcal D_k,F)$ be the disjoint union of the
graphs in $\hat{\mathcal D}_k\cap\Forb(F)$.
In this situation we can use  Theorem~\ref{locdual} and put $U(\mathcal K,\mathcal F)=U(\mathcal D_k,F)^{\Uparrow
  (p+1)\dots\Uparrow N}$.
\end{proof}  

\section{Concluding remarks}

{\bf 1. On Hadwiger conjecture}

Let us list the following corollary of Theorem \ref{main1}

\begin{corollary}
\label{trivi}
Let $\mathcal K$ be a proper minor closed class of graphs. 
Let $\mathcal F$ be a finite set of finite connected  graphs.
Then there exists a finite graph $U(\mathcal K,\mathcal F)\in\Forb(\mathcal
F)$ such that any graph of $\mathcal K\cap\Forb(\mathcal F)$ has
a homomorphism to $U(\mathcal K,\mathcal F)$.
\end{corollary}

It is well known (\cite{HN}) that one can reformulate the Hadwiger conjecture as the existence 
of a maximum (in the homomorphism order)  for  every proper minor closed class. 
Let $h = h(\mathcal K)$ be the Hadwiger number of the class $\mathcal K$. Then 
$K_{h+1} \not\in \mathcal K$ and Corollary \ref{trivi} gives at least a $K_{h+1}$-free bound
of the class $\mathcal K$. In fact we can get a bound with any set of the same local properties as the class $\mathcal K$ itself.

{\bf 2. On bounded expansion classes}

Let $\mathcal K$ be the class of all graphs $G$ which have bounded expansion 
with the expansion function $f$. Formally, 
$\mathcal K = \{G; \rdens{r}(G) \leq f(r), r = 1,2,\ldots\}$.
Assume that $p$ is minimal with $K_p \not\in \mathcal K$.
Then $\rdens{0}(G) \leq p-1$ for every $G \in \mathcal K$. Thus every $G \in  \mathcal K$ 
is $p-1$-degenerated,  If $f$ is monotonne then also $K_p \in \mathcal K$ and thus $\mathcal K$
has maximum. Thus Hadwiger conjecture holds for 
bounded expansion classes determined by a monotonne expansion function.

Note also that for constant expansion functions the bounded expansion classes are 
proper just  minor closed classes. This may be seen as follows:

Assume that $\mathcal K$ is bounded expansion class bounded by a constant function $C$.
Explicitely, we assume that for every $G \in \mathcal K$ holds  
$\rdens{r}(G) \leq C$ and thus also $\dens(G) \leq C$. Let $H$ be a minor of $G$. Then obviously
$\dens(H) \leq \dens(G) \leq C$ and thus  $\mathcal K$ is minor closed.
It follows that $K_{2C+1}$ is a forbidden minor of $\mathcal K$.

{\bf 3. On distal colorings - exact powers}

We now explain a particular consequence of our main result in a
greater detail. Let $G$ be a graph, $p$ a positive integer. Denote
by $G^{\natural p}$ the graph $(V, E^{\natural p})$ where $\{x, y\}$ is an edge of
$E^{\natural p}$ iff there exists a path $P$ in $G$ from $x$ to $y$ of length
$p$. The graph $G^{\natural p}$
 could be called
{\em exact $p$-power} of $G$. Clearly graphs $G^{\natural 2}$ and
all graphs $G^{\natural p}, p$ even, may  have  unbounded
chromatic number even for the case of trees (consider  subdivision
stars), and the only (obvious) bound is $\chi(G^{\natural p})\leq
\Delta(G)^p+1$. Similarly, for any odd $p$ there are $3$-colorable
graphs  $G$ for which is the chromatic number $\chi(G^{\natural
p})$ may be arbitrarily large. However for $p$ odd and arbitrary
proper minor closed class and even class with bounded expansion 
we have the following (perhaps
surprising):

\begin{theorem}
\label{ppower}
For any  class $\mathcal K$ with bounded expansion and for every  odd integer $p\geq 1$,
there exists an integer $N=N(\mathcal K,p)$ such that all the graphs $G^{\natural p}, G \in {\mathcal K}$ and
${\rm odd\!-\!girth}(G)>p$ have chromatic number $\leq N$:
For any $G\in\mathcal K$,
\begin{equation*}
{\rm odd\!-\!girth}(G)>p \Longrightarrow \chi(G^{\natural p})\leq N
\end{equation*}
% For every $p$ odd and every proper minor closed class $\mathcal K$ there exists an integer
% $N = N(p, {\mathcal K})$ such that all the graphs $G^{(p)}, G \in {\mathcal K}$ and
% ${\rm odd\!-\!girth}(G)>p$ have chromatic number $\leq N$.
\end{theorem}

Theorem \ref{ppower} follows immediately from Theorem \ref{main1}. It suffices to consider 
$\mathcal{F} = \{C_{p}\}$. In this case any graph $U(\mathcal K,\mathcal F)\in\Forb(\mathcal
F)$ and any homomorphism $c: G \xhom{} U(\mathcal K,\mathcal F)$ gives a desired coloring 
by $N = \card(V(U(\mathcal K,\mathcal F))$ colors.

With a little more care one can prove the following result which we state without proof 
(as it may be generalized in yet another direction, see \cite{POMNIV}. 
We take time for a definition of {\em exact distance } graph:
Let $G$ be a graph, $p$ a positive integer. Denote
by $G^{[\natural p]}$ the graph $(V, E^{[\natural p]})$ where $\{x, y\}$ is an edge of
$E^{[\natural p]}$ iff the distance of $x$ and $y$  in $G$ is $p$.

\begin{theorem}
\label{ppower1}
For any  class $\mathcal K$ with bounded expansion and for every  odd integer $p\geq 1$,
there exists an integer $N'=N(\mathcal K,p)$ such that all the graphs 
$G^{[\natural p]}, G \in {\mathcal K}$  have chromatic number $\leq N$.
\end{theorem}

Note that in both Theorems \ref{ppower},\ref{ppower1} we cannot replace the conditions in the definition of powers $G^{\natural p},G^{[\natural p]}$  by the existence of a path (or even induced 
path) of length $p$. See \cite{POMNIV} for a more detailed disscussion).
 
{\bf 4. On universality of posets.}
It follows from \cite{HubN} 
that the class $\mathcal SPG$ of all finite series paralel graphs is universal partial 
ordered class. What this means (using non-trivial result that the homomorphism order of all
finite graphs is universal, \cite{HubN,HN}) 
is that to every finite graph $G$ we can associate a series 
parallel graph $\Phi(G)$ such that for any two graphs $G, H$

$$ G \leq H \iff \Phi(G) \leq \Phi(H).$$

Thus bounded tree width graphs form a homomorphism universal class. Note that for the 
tree depth such a statement  is deeply not true as we cannot even find an infinite antichain.
In fact up to homomorphism equivalence the class of all graphs with a bounded tree depth
is a finite class (see  Lemma \ref{finitness}). Also this indicate that low tree-depth partitions
are much more restrictive than low tree width partitions.

{\bf 5. Regular partitions.}

Implicit in our proof of restricted dualities is the following partition 
result.
By Lemma \ref{finitness} there exists a finite set 
$\hat{\mathcal D}_p$ of graphs with tree depth  $p$, so that
any graph with tree-depth  $p$ is hom-equivalent to one graph
in the set. This implies: 
\begin{theorem}
\label{szemeredlike}
For every class ${\mathcal K}$ with bounded expansion and for every 
positive integer $p$ there exists a positive integer $N = N(\mathcal K, p)$
such that for every graph $G \in \mathcal K$ there exists a coloring 
$V(G) = V_1 \cup \ldots \cup V_N$ such that the
subgraph $G_J$  of $G$ induced by any $j = \card{J}$ classes has 
the following property: 
\begin{itemize}
\item each  component of $G_J$ is homomorphism aequivalent to one of the
graphs in the finite set $\hat{\mathcal D}_p$.
 \end{itemize}
\end{theorem}

This stronger decomposition theorem may be used for an alternative proof of 
Theorem \ref{main1}. Finally note that our results may be regarded as 
(very) regular partitions of graphs with bounded expansion. While the 
celebrated Szemeredi regularity lemma \cite{sz} applies to  dense graphs
the classes of bounded expansion are on the other side of spectrum: their
edge densities are (hereditarily) small. For these classes one can the achieve 
a very regular partitions, 
essentially describing all components  which may occur in graphs induced by a few color classes. 

%\bibliographystyle{amsplain}
%\bibliography{bib,biblio}

\providecommand{\bysame}{\leavevmode\hbox to3em{\hrulefill}\thinspace}
\providecommand{\MR}{\relax\ifhmode\unskip\space\fi MR }
\providecommand{\MRhref}[2]{%
  \href{http://www.ams.org/mathscinet-getitem?mr=#1}{#2}
}
\providecommand{\href}[2]{#2}

\end{document}